\documentclass[12pt, reqno]{amsart}
\usepackage{amsmath, amsthm, amscd, amsfonts, amssymb, graphicx, hyperref,color}
\newtheorem{theorem}{Theorem}[section]

\newtheorem{proposition}[theorem]{Proposition}
\newtheorem{corollary}[theorem]{Corollary}
\theoremstyle{definition}

\newtheorem{example}[theorem]{Example}

\theoremstyle{remark}
\newtheorem{remark}[theorem]{Remark}
\numberwithin{equation}{section}
\begin{document}
%%%%%%%%%%%%%%%%%%%%%%%%%%%%%%%%%%%%%%%%%%%%%%%%%%
%%%%%%%%%%%%%%%%%%%%%%%%%%%%%%%%%%%%%%%%%%%%%%%%%
\setcounter{page}{1}
%%%%%%%%%%%%%%%%%%%%%%%%%%%%%%%%%%%%%%%%%%%%%

%-------------------------- Pleased do not change the following line-------------------------------------------
%\noindent \textcolor[rgb]{0.99,0.00,0.00}{}\\[.5in]
%--------------------------------------------------------------------------------------------------------------

\title[A factorization result for quadratic pencils of accretive operators ]
      {A factorization result for quadratic pencils of accretive operators and applications}
       
%%%%%%%%%%%%%%%%%%%%%%%%%%%%%%%%%%%%%%%%%%%%%%%%%
%%%%%%%%%%%%%%%%%%%%%%%%%%%%%%%%%%%%%%%%%%%%%%%%%
\author[ F. Bouchelaghem, M. Benharrat]{Fairouz Bouchelaghem$^1$, Mohammed Benharrat$^{2*}$}
%%%%%%%%%%%%%%%%%%%%%%%%%%%%%%%%%%%%%%%%%%%%%%%
%%%%%%%%%%%%%%%%%%%%%%%%%%%%%%%%%%%%%%%%%%%%%%%%
\address{$^{1}$ 
D\'{e}partement de Math\'{e}matiques, Univérsit\'{e} Oran1-Ahmed Ben Bella, BP 1524 Oran-El M'naouar, 31000 Oran, Alg\'{e}rie.}
\email{\textcolor[rgb]{0.00,0.00,0.84}{fairouzbouchelaghem@yahoo.fr}}
%%%%%%%%%%%%%%%%%%%%%%%%%%%%%%%%%%%%%%%%%%%%%%%%%%%%%
\address{$^{2}$ D\'{e}partement de G\'{e}nie des syst\'{e}mes,
	Ecole Nationale Polytechnique d'Oran-Maurice Audin (Ex. ENSET d'Oran), 
	BP 1523 Oran-El M'naouar, 31000 Oran, Alg\'{e}rie.}
\email{\textcolor[rgb]{0.00,0.00,0.84}{	mohammed.benharrat@enp-oran.dz, mohammed.benharrat@gmail.com}}

%%%%%%%%%%%%%%%%%%%%%%%%%%%%%%%%%%%%%%%%%%%%%%%%%%%%%%%%%%%%

%\dedicatory{This paper is dedicated to Professor ABCD}



\date{01/10/2020.
	\newline \indent $^{*}$ Corresponding author\\
	This work was supported by the Laboratory of Fundamental and Applicable Mathematics of Oran (LMFAO) and the Algerian research project: PRFU, no. C00L03ES310120180002.}
%%%%%%%%%%%%%%%%%%%%%%%%%%%%%%%%%%%%%%%%%%%%%%%%%%%%%%%%%%%%%%%%
\begin{abstract}
%%%%%%%%%%%%%%%%%%%%%%%%%%%%%%%%%%%%%%%%%%%%%%%%%%%%%%%%%%%%%
A canonical factorization is given for  a quadratic  pencil of accretive operators in a Hilbert space. We  establish a criterion in order that the linear factors, into which the pencil splits,  generates an holomorphic semi-group of contraction operators.  As an application, we study a result of existence, uniqueness, and maximal regularity of the strict solution  for complete abstract second order differential equation  in the non-homogeneous case. An illustrative example is also given.
\end{abstract}
%%%%%%%%%%%%%%%%%%%%%%%%%%%%%%%%%%%%%%%%
%%%%%%%%%%%%%%%%%%%%%%%%%%%%%%%%%%%%%%%%%%%%%%%%%%%%%%%%%%%%%%%
 \maketitle
%%%%%%%%%%%%%%%%%%%%%%%%%%%%%%%%%%%%%%%%%%%%%%%%%%%%%%%%%%%%%%
%%%%%%%%%%%%%%%%%%%%%%%%%%%%%%%%%%%%%%%%%%%
%Introduction
%%%%%%%%%%%%%%%%%%%%%%%%%%%%%%%%%%%%%%%%%%

\section{introduction}
%%%%%%%%%%%%%%%%%%%%%%%%%%%%%%%%%%
Many problems in   mathematical physics and mechanics can be described  by the following second order  linear differential equation
%In this paper we consider a linear second order differential equation
\begin{equation}\label{MainEquation}
	u''(t)-2Bu'(t)-Cu(t)=0, %\quad t\geq 0
\end{equation}
where $u(t)$ is a vector-valued function in an appropriate (finite or infinite dimensional) Hilbert space $\mathcal{H}$, $B$ and $C$
are linear (bounded or unbounded) operators on $\mathcal{H}$. Properties of the differential equation \eqref{MainEquation} are
closely connected with spectral properties of a quadratic pencil
\begin{equation}\label{qop0}%%%%%%%%%%%%%%%%%%%%%%%
Q(\lambda)=\lambda^2 I-2 \lambda B-C, \qquad (\lambda \in \mathbb{C});
\end{equation}
%%%%%%%%%%%%%%%%%%
which is obtained by substituting exponential functions $u(t)=\exp(\lambda t)x$, $x\in \mathcal{H}$ into \eqref{MainEquation}.
In many applications $B$ and $C$ are  self-adjoint positive definite operators. An  important and subtle problem in the theory of such operator pencils is to factoring them and studying the spectral properties of the factors. Under some  assumptions, Krein and  Langer \cite{KreinL64} proved that a self-adjoint polynomial of the form \eqref{qop0} can always be written  as a product of two linear factors as follows
\begin{equation}\label{fact} \lambda^2 I-2 \lambda B-C =(\lambda I- Z_1)(\lambda I-Z_2),
\end{equation}
with   $Z_1$ and $Z_2$  are a roots of the quadratic operator equation
\begin{equation}\label{qope}
	Q(Z)=Z^2 I-2  BZ-C=0.
\end{equation}
 Of particular interest is the separation  of spectral values of $Q$ between the  spectra of the  roots. Such separation may be complicated, 
 even in the case of eigenvalues, see   \cite{Shkalikov89} and references therein.  The
 factorization theorems  have  been studied extensively also for the self-adjoint quadratic operator pencils under the extra condition of strong and weak  damping. For the exhaustive survey on these topics, please see the two seminal books \cite{Markus88} and \cite{ Moller2015} and the references therein.

 But some models of continuous mechanics are reduced to differential equation \eqref{MainEquation} with sectorial
 operators, see \cite{AglazinKiiko1, Artamonov1, Favini2008, IlyushinKiiko1} and references therein. In this cases methods,
 developed for self-adjoint operators, cannot be applied.

The main objective of the manuscript is to find sufficient conditions on, in general, unbounded linear accretive operators $B$ and $C$ on the Hilbert space, under which a factorization \eqref{fact} is possible.
The approach is  based on the perturbation theory of  accretive operators. We also obtain a
criterion in order that the linear factors, into which the pencil splits,  generates an holomorphic semi-group of contraction operators. We apply this result to establish a theorem  of existence, uniqueness, and maximal regularity of the strict solution  of an abstract second order evolutionary equations generated by such pencils in the non-homogeneous case.

%%%%%%%%%%%%%%%%%%%%%%%%%%%%%%%%%%%%%%%%%%%%%%%%%%%%%%%%%%%%%
%%%%%%%%%%%%%%%%%%%%%%%%%%%%%%%%%%%%%%%%%%%%%%%%%%%%%%%%
\section{Accretive operators framework} 
%%%%%%%%%%%%%%%%%%%%%%%%%%%%%%%%%%%%%%%%%%%%%%%%%%%%%%%%%%%
In this section we introduce the notation and the operator theoretic framework used in the rest of our work. Throughout this paper  $\mathcal{H}$ is a complex   Hilbert space with inner product $<\cdot, \cdot>$ and norm $\|\cdot\|$. Let $\mathcal{B} (\mathcal{H})$  denote the Banach space of all bounded linear operators on $\mathcal{H}$. Given a linear operator $T$ on $\mathcal{H}$ we denote by  $\mathcal{D} (T)$, $\mathcal{N}(T)$, and $\mathcal{R}(T)$ the domain, the null space and the range of $T$, respectively. For a closable densely defined linear
operator $T$ in some Hilbert space $\mathcal{H}$ we denote by $\rho (T)$ the resolvent set, by  $\sigma (T)=\mathbb{C}\backslash \rho (T)$ the spectrum, and by $\sigma_p(T)$ the point spectrum of $T$.
 For $\lambda \in \rho (T),$ the inverse $(\lambda I-T)^{-1}$ is, by the closed graph theorem, a bounded operator on $\mathcal{H}$ and will be called the
resolvent of $T$ at the point $\lambda$.   
%%%%%%%%%%%%%%%%%%%%%%%%%%%%%%%%%%%%%%%%%%%%%%%%%%%%%%%%%%%%%%%%%%%%%%%%%%%

Recall that  a linear operator $T$ with domain $\mathcal{D}(T)$ in a complex Hilbert space $\mathcal{H}$ is said to be accretive if
$$
Re<Tx,x>\geq 0 \qquad \text{ for all } x\in \mathcal{D}(T)
$$
or, equivalently if
$$\Vert(\lambda +T)x\Vert\geq \lambda \Vert x\Vert \qquad \text{ for all } x\in \mathcal{D}(T) \text{ and  }\lambda>0.$$
An accretive operator $T$ is called \textit{maximal accretive}, or
\textit{$m$-accretive} for short, if one of the following equivalent conditions is satisfied:
\begin{enumerate}
	\item $T$ has no proper accretive extensions in $\mathcal{H}$;
	\item $T$ is densely defined and $\mathcal{R}(\lambda +T)=\mathcal{H}$ for some (and hence for every) $\lambda>0$;
	\item $T$ is densely defined and closed, and $T^*$ is accretive.
\end{enumerate}
 In particular, every m-accretive operator is accretive
and  closed densely defined, its adjoint is also m-accretive (cf. \cite{Kato}, p. 279). Furthermore, $$ (\lambda + T)^{-1} \in \mathcal{B} \ (\mathcal{H}) \quad  \text{ and } \quad \left\|(\lambda + T)^{-1}\right\| \leq \frac{1}{\lambda}   \text{ for } \lambda >0.$$
In particular, a bounded accretive  operator  is m-accretive.  \^Ota showed in  \cite[Theorem 2.1]{Ota1984} that, if $T$ is closed and  an accretive such that there is a positive integer $n$ with $\mathcal{D}(T^n)$ is dense in $\mathcal{H}$ and $\mathcal{R}(T^n)\subset \mathcal{D}(T)$, then $T$ is bounded . In particular, for a  closed   and  accretive operator $T$,   if $\mathcal{R}(T)$ is contained in $\mathcal{D}(T)$, or in  $\mathcal{D}(T^*)$, then $T$ is automatically bounded, see also \cite[Theorem 3.3]{Ota1984}. Also, if $T$ is maximal accretive, then
\begin{equation}
\label{kernel-accretive}
\mathcal{N} (T) = \mathcal{N} (T^*) \quad  \text{ and } \quad \mathcal{N} (T)\subseteq\mathcal{D} (T) \cap \mathcal{D} (T^*).
\end{equation} 
%%%%%%%%%%%%%%%%%%%%%%%%%%%%%%%%%%%%%%%%%%%%%%%%%%%%%%%

The numerical range is very useful set by what we can  we define the  accretive  operators. For a linear operator  $T: \mathcal{D}  (T)\rightarrow  \mathcal{H}$  it is defined by
\begin{equation}
W(T):=\{ <Tx, x> : \quad  x\in \mathcal{D}  (T) , \quad \text{ with } \left\|x\right\|=1\},
\end{equation}
It is well-known that  $W(T)$ is a convex set of the complex plane (the Toeplitz-Hausdorff
theorem), and in general is neither open nor closed, even for a closed operator $T$. Clearly,  an operator $T$  is accretive when   $W(T)$ is contained in the closed right half-plane
\[
W(T) \subset \overline{\mathbb{C}_+} := \{ z\in\mathbb{C} : \rm Re (z) \ge 0\}.
\]
 Further,  if $T$ is m-accretive operator then $W(T)$ has the so-called spectral inclusion property
\begin{equation}\label{specinclusion}
\sigma (T)\subset\overline{ W(T)}.
\end{equation}
%%%%%%%%%%%%%%%%%%%%%%%%%%%%%%%%%%%%%%%%%%%%%%%%%%%%%%%%%%

%A linear operator $T$ in a Hilbert space $H$ is called \textit{dissipative}, or \textit{$m$-dissipative} for short, if
%$-I -T$ is accretive or $m$-accretive, respectively; in this case the numerical range of $T$ is contained in
%the closed upper half-plane, $W(T)\subset \{ z\in \mathbb{C}: Im z \ge 0\}$.

%An accretive operator $T$ is called \emph{coercive} if there exists $m>0$ with $Re <Tx, x> \ge m\|f\|^2$ for all $x\in \mathcal{D}  (T)$.

A linear operator $T$ in a Hilbert space $H$ is called \emph{sectorial} with vertex $z=0$ and
semi-angle $\omega \in [0,\pi/2)$, or \emph{$\omega$-accretive} for short, if its numerical range is contained in a
closed sector with semi-angle~$\omega$,
\begin{equation}
\label{sectorial}
W(T) \subset \overline{\mathcal{S}(\omega)}:=\left\{z\in\mathbb{C} : |\arg z|\leq \omega \right\}
\end{equation}
or, equivalently, $$|Im<Tx, x> |\!\le\! \tan\omega \, Re<Tx, x> \qquad \text{ for all }  x \in  \mathcal{D}  (T).$$

An $\omega$-accretive  operator  $T$ is called   m-$\omega$-accretive, if it is $m$-accretive. We have $T$ is m-$\omega$-accretive if and only if the operators $e^{\pm i\theta} T$ is  m-accretive for $\theta=\frac{\pi}{2}-\omega$, $0 < \omega\leq  \pi/2$.
The resolvent set of an m-$\omega$-accretive operator $T$ contains the
set $\mathbb{C} \setminus \overline{\mathcal{S}(\omega)}$ and
\[
\|(T-\lambda I)^{-1}\| \le \cfrac{1}{{\rm dist}\left(\lambda,\mathcal{S}(\omega)\right)},
\quad \lambda\in \mathbb{C}\setminus \overline{\mathcal{S}(\omega)}.
\]
%%%%%%%%%%%%%%%%%%%%%%%%%%%%%%%%%%%%%%%%%%%%%%
In particular,  m-$\pi/2$-accretivity means m-accretivity.  A $0$-accretive operator
is symmetric. An operator is positive if and only if it is m-$0$-accretive. 
%%%%%%%%%%%%%%%%%%%%%%%%%%%%%%%%%%%%%%%%%%%%%%%

It is known that the $C_0$-semigroup $\mathcal{T}(t)=\exp(-tT)$, $t\ge 0$, has contractive  and holomorphic continuation into the sector $ \overline{\mathcal{S}(\pi/2\omega)}$ if and only if the generator $T$  is m-$\omega$-accretive, see \cite[Theorem V-3.35] {Kato}.

%%%%%%%%%%%%%%%%%%%%%%%%%%%%%%%%%%%%%%%%%%%%%%%%%%%%%%%
Recently, the authors of \cite{Arlinskii2010}  obtained a precise localization of the numerical range of one-parameter semigroup $\mathcal{T}(t)=\exp(-tT)$, $t\ge 0$, generated  by   an  m-$\omega$-accretive operator, $\omega \in [0,\pi/2)$. More precisely, by \cite[Theorem 3.4]{Arlinskii2010}, we have
\begin{equation}
W(\exp(-tT))\subseteq \Omega(\omega)=\{z\in \mathbb{C}: \left|Im \sqrt{z} \right|\leq \dfrac{1}{2}(1-\left| z\right| ) \tan(\omega)\}, \quad t\geq 0,
\end{equation}
with limiting cases:  $\Omega(0)=[0,1]$ and $\Omega(\pi/2)=\overline{\mathbb{D}}$. In particular, the family ${\exp(-tT)}_{t\geq 0}$ is a quasi-sectorial contractions semigroup in the terminology of \cite{Arlinskii2010}.
%%%%%%%%%%%%%%%%%%%%%%%%%%%%%%%%%%%%%%%%%%%%%%%%%%%%%%%

We mention that if $T$ is $m$-accretive, then for each $\alpha\!\in\!(0,1)$  the fractional powers $T^\alpha$, $0<\alpha<1$,  are defined by the following Balakrishnan formula, see \cite{Balakrishnan60},
$$
T^\alpha x =\dfrac{\sin (\pi \alpha )}{\pi} \int_{0}^{\infty}\lambda^{\alpha -1}T(\lambda +T)^{ -1}xdt,
$$
for all $x\in \mathcal{D}  (T)$. The operators $T^{\alpha}$ are m-$(\alpha\pi)/2$-accretive  and, if $\alpha\!\in\! (0,1/2)$, then
$\mathcal{D} (T^{\alpha})= \mathcal{D} (T^{*\alpha})$.  It was proved in \cite[Theorem 5.1]{Kato61} that, if $T $ is $m$-accretive, then
$\mathcal{D} (T^{1/2}) \cap \mathcal{D}(T^{*1/2})$ is a core of both $T^{1/2}$  and $ T^{*1/2}$ and the real part $Re T^{1/2} :=\!(T^{1/2}\!+\!T^{*1/2})/2 $ defined on $\mathcal{D} (T^{1/2})\cap  \mathcal{D} (T^{*1/2})$ is a selfadjoint operator.
Further,~by~\cite[Corollary ~2]{Kato61},
\begin{equation}
\label{ravhalf}
\mathcal{D}(T)= \mathcal{D}(T^*) \implies \mathcal{D} (T^{1/2})=\mathcal{D} (T^{*1/2})=\mathcal{D} (T^{1/2}_R) = \mathcal{D}[\phi],
\end{equation}
where $\phi$ is the closed form associated with the sectorial operator $T$ via the first representation theorem \cite[Sect.~VI.2.1]{Kato}
and $T_R$ is the non-negative selfadjoint operator associated with the real part of $\phi$ given by
${\rm Re}\, {\phi} := (\phi + \phi^*)/2$.
%%%%%%%%%%%%%%%%%%%%%%%%%%%%%%%%%%%%%%%%%%%%%%%%%%%%%%%%%

%%%%%%%%%%%%%%%%%%%%%%%%%%%%%%%%%%%%%%%%%%%%%%%%%%%%%%%%%%%%%%
\section{A  canonical factorization of a monic quadratic operator pencils}\label{sec:R} 
%%%%%%%%%%%%%%%%%%%%%%%%%%%%%%%%%%%%%%%%%%%%%%%%%%%%%%%%%%%%%%
In this section we will investigate a canonical factorization of quadratic operator pencils $\mathcal{Q}$ of
the form
\begin{equation}\label{qop}
\mathcal{Q}(\lambda)=\lambda^2 I-2 \lambda B-C,
\end{equation}
on a Hilbert space with domain $\mathcal{D} (\mathcal{Q})=\mathcal{D} (B) \cap \mathcal{D} (C)$, where $\lambda\in \mathbb{C}$ is the spectral parameter and the two operators  $B$ and $ C $   with domain $\mathcal{D} (C)$ and $\mathcal{D} (B)$, respectively, satisfy one of the following conditions,

\begin{itemize}
	\item [(\textbf{C.1})] there exists $\alpha\geq 0$, $0\leq \beta <1$ and $\delta\geq 0 $ such that
	\begin{equation*}\label{bmaccrcb}
	Re<B^2x, Cx> \geq-\alpha\left\|x\right\|^2- \beta\left\|B^2x\right\|^2-\delta\left\|B^2x\right\|\left\|x\right\|,
	\end{equation*}
	for all   $x\in \mathcal{D} (B^2)\subset \mathcal{D} (C)$.
	\item [(\textbf{C.2})]  $C$ is $B^2$-bounded with lower bound $< 1$. i.e. there exists $a\geq 0$ and $0\leq b <1$   such that
	\begin{equation*}\label{Tbounded}
	\left\|Cx\right\|^2\leq a \left\|x\right\|^2+b\left\|B^2x\right\|^2,\quad  \text{ for all }  x\in \mathcal{D} (B^2)\subset \mathcal{D} (C).
	\end{equation*}
	\item [(\textbf{C.3})] $I+C(B^{2}+t_0)^{-1}$ is boundedly invertible, for some $t_0>0$.
	\item [(\textbf{C.4})] $B$ is accretive and $\mathcal{D} (B)\subset \mathcal{D} (C)$.
	\item [(\textbf{C.5})] $B$ is accretive and $C$ is bounded.
\end{itemize}

%%%%%%%%%%%%%%%%%%%%%%%%%%%%%%%%%%%%%%%%%%%%%%%%%%%%%%%%%%%%%
\begin{proposition}\label{thm:1D1} Let $B^2$  be m-accretive and $C$  is accretive. If	the operator $B$ and $C$ verifies one of the conditions above, then the operator  $\Lambda=B^2 +C$ with domain $ \mathcal{D}  (B^2)$  is  m-accretive.	
\end{proposition}
%%%%%%%%%%%%%%%%%%%%%%%%%%%%%%%%%%%%%%%%%%%%%%%%%%%%%%%%%%%%%

\begin{proof}Assume  (\textbf{C.1}), then by  \cite[Theorem 3.10]{Okazawa82} we can prove that $B^2 +C$ is m-accretive. For the convince of the reader we give a detailed proof of this fact and adapted to the Hilbert case. First, we have $B^2+C$ is accretive densely defined. We   show that $B^2+C$  is closed. In fact, it follows from (\textbf{C.1}) that
	\begin{align*}
	\left\|B^2x\right\|^2& = Re<B^2x, B^2x>\\
	&\leq  Re<(B^2+C)x, B^2x>+ \alpha\left\|x\right\|^2+ \beta\left\|B^2x\right\|^2+\delta\left\|B^2x\right\|\left\|x\right\|,
	\end{align*}
for all   $x\in \mathcal{D} (B^2)$.	So, we have
	$$(1-\beta) \left\|B^2x\right\|^2\leq [\delta\left\|x\right\|+\left\|(B^2+C)x\right\|] \left\|B^2x\right\|+ \alpha\left\|x\right\|^2,$$
	for all   $x\in \mathcal{D} (B^2)$. Solving this inequality, we obtain
	\begin{equation}\label{1-beta}
	\left\|B^2x\right\|\leq \dfrac{1}{1-\beta}\left\|(B^2+C)x\right\|+\kappa \left\|x\right\| 
	\end{equation}
	for all   $x\in \mathcal{D} (B^2)$, with $ \kappa = \dfrac{\alpha +\sqrt{\delta(1-\beta)}}{1-\beta}$. On the other hand, since $\mathcal{D} (B^2)\subset \mathcal{D} (C)$, with $\mathcal{D} (B^2)$ dense in $\mathcal{H}$, there exists a constant $\vartheta>0$, such that 
	\begin{equation}\label{closable}\left\|Cx\right\|\leq \vartheta (\left\|x\right\|+\left\|B^2x\right\|),	\end{equation}
	for all $ x\in \mathcal{D}  (B^2)$. Now, let a sequence $ (x_n)_n\subset\mathcal{D} (B^2) $ such that $x_n\longrightarrow x $ and $(B^2+C)x_n\longrightarrow y$. Applying the inequality \eqref{1-beta} to $x$ replaced by $x_n-x_m$, we see that the sequence $(B^2x_n)_n$ converge. Since $B^2$ is closed we conclude that $B^2x_n\longrightarrow B^2x$ and $x\in \mathcal{D} (B^2)$. By \eqref{closable},  we have
	\begin{equation*}\left\|C(x_n-x)\right\|\leq \vartheta \left\|x_n-x\right\|+(1+\vartheta)\left\|B^2(x_n-x)\right\|.	\end{equation*}
	Hence $(B^2+C)x_n\longrightarrow(B^2+C)x$ and $y=(B^2+C)x$, which shows $B^2+C$ is closed.  On the other hand, we have
	\begin{equation*}\label{bmaccrcb}
	Re<B^2x, tCx> \geq-t\alpha\left\|x\right\|^2- t\beta\left\|B^2x\right\|^2-t\delta\left\|B^2x\right\|\left\|x\right\|,
	\end{equation*}
	for all $0\leq t\leq 1$. Since $t\beta<1$,
	 by the same argument, we assert that $(B^2+tC)$ is closed for all $0\leq t\leq 1$. Hence, $(B^2+tC)$ is closed and accretive for all $0\leq t\leq 1$. By  \cite[Lemma  3.1]{Okazawa82},  the dimension of $\mathcal{R}(B^2+C+\lambda)^{\bot} $ and $\mathcal{R}(B^2+\lambda)^{\bot}$ are the same for all $\lambda>0$. Since $B^2$ is m-accretive, we conclude that $B^2+C$ is m-accretive.

	If (\textbf{C.2}) holds, the result follows by \cite[Theorem 2.]{Gustafson66}. Also (\textbf{C.2}) implies (\textbf{C.1}) in the case of $\gamma =0$, see \cite[Remark 4.4]{Okazawa73}. In fact, setting $\alpha=a/2$ and $\beta=(b+1)/2$, we have that for $x\in \mathcal{D} (B^2)$,
	\begin{align*}
	\left\|Cx\right\|^2&\leq 2\alpha \left\|x\right\|^2+(2\beta-1)\left\|B^2x\right\|^2\\
	&\leq \left\|(B^2+C)x\right\|^2-\left\|B^2x\right\|^2+2\alpha \left\|x\right\|^2+2\beta\left\|B^2x\right\|^2\\
	&= 2(Re<B^2x, Cx>+ \alpha\left\|x\right\|^2+ \beta\left\|B^2x\right\|^2)+\left\|Cx\right\|^2.
	\end{align*}
	Hence $Re<B^2x, Cx>+ \alpha\left\|x\right\|^2+ \beta\left\|B^2x\right\|^2\geq 0$.
	
	Assume that (\textbf{C.3}) is satisfied. Since $B^2 +C$ is densely defined and accretive, it suffices to show that $\mathcal{R}(B^2+C+t_0) =\mathcal{H}$. But this follows immediately from 
	$$  B^2+C+t_0=(I+C(B^{2}+t_0)^{-1})(B^{2}+t_0),$$
	and clearly $B^2+C+t_0$ is invertible.
	
	Now, we consider (\textbf{C.4}). Since $B$ is an accretive operator, by \cite[Theorem 1.2]{Hayashi2017}, we have for an arbitrary $\nu>0$, 
	\begin{equation}\label{BBB2}
	\left\|B x\right\|^2 \leq  \nu \left\|x\right\|^2+\dfrac{1}{\nu}\left\|B^2 x\right\|^2,
	\end{equation}
	for all  $x\in  \mathcal{D}  (B^2)$. Since $\mathcal{D} (B)\subset \mathcal{D} (C)$, with $\mathcal{D} (B)$ dense in $\mathcal{H}$, there exists a constant $\eta>0$, such that 
	$$\left\|Cx\right\|^2\leq \eta \left\|Bx\right\|^2,$$
	for all $ x\in \mathcal{D}  (B)$. It follows that
	$$
	\left\|Cx\right\|^2\leq \eta (\nu\left\|x\right\|^2+\dfrac{1}{\nu}\left\|B^2 x\right\|^2),
	$$
	for all $ x\in \mathcal{D}  (B^2 )$. Choosing $\nu > 0 $ so large that $\dfrac{\eta}{\nu} < 1$, we get $C$ is $B^2$-bounded with lower bound $< 1$.	
	
	Finally, clearly (\textbf{C.5}) is  a particular case of (\textbf{C.4}).
\end{proof}
%%%%%%%%%%%%%%%%%%%%%%%%%%%%%%%%%%%%%%%%%%%%%%%%%%%%%%%%%%
\begin{remark}\label{remark3.2} \begin{enumerate}
		\item In  (\textbf{C.3}), if we assume further $\mathcal{D} (B^2)\subset \mathcal{D} (C)$, then by \cite[Proposition 2.12]{Yoshikawa72}, the lower bound $b$ in (\textbf{C.2}) is equal to $\sup_{t>0}\left\|C(B^{2}+t)^{-1} \right\|$. Hence, if we assume further, $\left\|C(B^{2}+t_0)^{-1} \right\|<1$  for some $t_0>0$, so $I+C(B^{2}+t_0)^{-1}$ is boundedly invertible, for some $t_0>0$. In this case (\textbf{C.3}) implies (\textbf{C.2}).
		\item 	If the condition (\textbf{C.4}) is satisfied, clearly $\mathcal{D} (\mathcal{Q})=\mathcal{D} (B)$.
		\item In (\textbf{C.4}) and (\textbf{C.5}), $B$ is m-accretive. Indeed, choosing $\nu > 0 $ so large that $\dfrac{1}{\nu} < 1$ in \eqref{BBB2}, we obtain $B$ is $B^2$-bounded with lower bound $< 1$. Then $B^2+B$ with domain $\mathcal{D}  (B^2 )$ is 	m-accretive. Now, let us remark that
		$$(\dfrac{1}{4}I+B^2+B)x=(\dfrac{1}{2}I+B)^2x$$
		for all  $x\in  \mathcal{D}  (B^2)$. Since the operator on the left-hand side is invertible, then $(\frac{1}{2}I+ B)^2$ is invertible, so  $\frac{1}{2}I+B$ is also invertible. It follows that $B$ is m-accretive.
	\end{enumerate}
\end{remark}
%%%%%%%%%%%%%%%%%%%%%%%%%%%%%%%%%%%%%%%%%%%%%%%%%%%%%%%%
\textbf{In the sequel, we assume that $B^2$  be m-accretive and $C$  is accretive verify the condition (\textbf{C.1}), unless otherwise specified. }

Now, we state some properties of the operator $\Lambda=B^2 +C$. The first important one is the existence and uniqueness of its  square root. This is an  immediate consequence of \cite[Theorem 3.35, p. 281]{Kato}.
%%%%%%%%%%%%%%%%%%%%%%%%%%%%%%%%%%%%%%%%%%%%%%%%%%%%%%%%%%%

%%%%%%%%%%%%%%%%%%%%%%%%%%%%%%%%%%%%%%%%%%%%%%%%%%%%%%%%%%
\begin{corollary}\label{root}
	The operator $\Lambda$ admits  unique square root $\Lambda^{\frac{1}{2}}$ m-$(\pi/4)$-accretive operator  with $ \mathcal{D}  (B^2)$  is a core of $\Lambda^{\frac{1}{2}}$ (that is, the closure of the restriction of  $\Lambda^{\frac{1}{2}}$ to  $ \mathcal{D}  (B^2)$ is again $\Lambda^{\frac{1}{2}}$).
\end{corollary}
%%%%%%%%%%%%%%%%%%%%%%%%%%%%%%%%%%%%%%%%%%%%%%%%%%%%%%%%%%%%

%%%%%%%%%%%%%%%%%%%%%%%%%%%%%%%%%%%%%%%%%%%%%%%%%%%%%%%%%%%%%%%%
\begin{proposition}If $C$ is $ \theta$-accretive, with $0\leq \theta <\pi/2$, then
	$$\mathcal{N}(\Lambda) \subset \mathcal{N}(B^2)\cap \mathcal{N}(C^*).$$
\end{proposition}
%%%%%%%%%%%%%%%%%%%%%%%%%%%%%%%%%%%%%%%%%%%%%Since $B^2$ and $\Lambda$ are   m-accretive operators, we know that  $\mathcal{N}(B^2)= \mathcal{N}(B) $ and $ \mathcal{N} (\Lambda)=\mathcal{N}(\Lambda^{\frac{1}{2}})$.
\begin{proof}
	(1)  Let $ x\in \mathcal{D}  (B^2)$, $x\neq 0$, such that $\Lambda x=0$, as before,  we have
	$$Re< \Lambda x,x>=Re< B^2x,x>+Re <C x,x>,$$
	then 
	$$Re< B^2x,x> \leq Re< \Lambda x,x> \quad \text{ and } \quad Re <C x,x> \leq Re< \Lambda x,x>.$$
	Therefore, $Re< \Lambda x,x>=0$ implies that  $Re<B^2 x,x>=0$ and  $Re< Cx,x>=0$. 
	On the other hand, since $C$ is $ \theta$-accretive, with $0\leq \theta <\pi/2$, then 
	$$\left|Im <Cx,x>\right| \leq tan(\theta) Re< Cx,x> .$$
	Thus,
	$$ Im <Cx,x>=0 \quad \text{ and } \quad  Im <B^{2}x,x>=-Im <Cx,x>=0,$$
	hence
	$$   <B^2x,x>=0 \quad \text{ and } \quad   <Cx,x>=0 .$$
	Since $B^2$ is  m-accretive operator, we conclude that $x\in \mathcal{N} (B^2)$ and $x \in \mathcal{N} (C^*)$.
\end{proof}

Now, we define the linear factors $Z_1$ and $Z_2$  by
$$Z_1=B+\Lambda^{\frac{1}{2}}$$
and
$$Z_2=B-\Lambda^{\frac{1}{2}}$$
with domain $\mathcal{D} (B)\cap \mathcal{D}(\Lambda^{\frac{1}{2}})$, into which the quadratic pencil \eqref{qop} can  be decomposed.

\begin{proposition}\label{thm:1D1F} Assume that $B(\mathcal{D}  (B^2))\subset \mathcal{D}  (B^2) $ and $\Lambda^{\frac{1}{2}}(\mathcal{D}  (B^2))\subset \mathcal{D}  (B^2) $. 
	$Q$ takes the following form,
	\begin{equation}\label{z1z2nc}
		Q(\lambda)x=\dfrac{1}{2}(\lambda I-Z_1)(\lambda I-Z_2)x+\dfrac{1}{2} (\lambda I-Z_2)(\lambda I-Z_1)x,
	\end{equation}
	for all  $x\in  \mathcal{D}  (B^2)$.
	
	In particular, if $B\Lambda^{\frac{1}{2}}=\Lambda^{\frac{1}{2}}B$ on $\mathcal{D} (B^2)$, then	$Q$ admits the following canonical factorization,	
	\begin{equation}\label{z1z2}
		Q(\lambda)x=(\lambda I-Z_1)(\lambda I-Z_2)x=(\lambda I-Z_2)(\lambda I-Z_1)x,
	\end{equation}
	for all  $x\in \mathcal{D}  (B^2)$.
	
\end{proposition}
%%%%%%%%%%%%%%%%%%%%%%%%%%%%%%%%%%%%%%%%%%%%%%%%%%%%%%%%%%%%%

\begin{proof} We have $\mathcal{D} (B^2)\subset \mathcal{D} (Z_1)=\mathcal{D} (Z_2)=\mathcal{D} (\Lambda^{\frac{1}{2}})\cap \mathcal{D} (B)$ and $\mathcal{D} (B^2)\subset \mathcal{D} (C)$. 
	
	The fact that $B(\mathcal{D} (B^2))\subset \mathcal{D} (B^2)$ and $\Lambda^{\frac{1}{2}}(\mathcal{D} (B^2))\subset \mathcal{D} (B^2)$, we have $\mathcal{D} (B^2)\subset \mathcal{D}(B\Lambda^{\frac{1}{2}})$, $\mathcal{D} (B^2)\subset \mathcal{D}(\Lambda^{\frac{1}{2}}B)$ and
	$\mathcal{D} (B^2)\subset \mathcal{D}(Z_1^2)$. Now,  we can  verify  that
	\begin{equation*}\label{qopeq}
		Z_1^2x -B Z_1x-Z_1Bx-Cx=0,
	\end{equation*}
	for all  $x\in  \mathcal{D}  (B^2)$, hence on $ \mathcal{D}  (B^2)$, we have
	\begin{align*}
		Q(\lambda ) &=Q(\lambda)-(Z_1^2 -B Z_1-Z_1B-C)\\
		&= \lambda^2 I-2 \lambda B-C - Z_1^2 +B Z_1+Z_1B+C\\
		& = \lambda^2 I - Z_1^2 -B(\lambda -Z_1 )- (\lambda -Z_1 )B\\
		&=\dfrac{1}{2}(\lambda -Z_1 )(\lambda I + Z_1-2B)+\dfrac{1}{2}(\lambda I + Z_1-2B)(\lambda -Z_1 )\\
		&=\dfrac{1}{2}(\lambda I-Z_1)(\lambda I-Z_2)+\dfrac{1}{2} (\lambda I-Z_2)(\lambda I-Z_1).
	\end{align*}
	Now,  if $B\Lambda^{\frac{1}{2}}=\Lambda^{\frac{1}{2}}B$ on $\mathcal{D} (B^2)$, we  obtain \eqref{z1z2}.	
\end{proof}
%%%%%%%%%%%%%%%%%%%%%%%%%%%%%%%%%%%%%%%%%%%%%%%%%%%%%%%%

In the sequel, we investigate some properties of the operators    $Z_1$ and $Z_2$.

%%%%%%%%%%%%%%%%%%%%%%%%%%%%%%%%%%%%%%%%%%%%%%%%%%%%%%%%%%%%%%%
\begin{proposition}\label{thm:1D1S}  Assume  that  $\mathcal{D}(\Lambda^{\frac{1}{2}}) \subset \mathcal{D} (B) $, then  for any $\varepsilon>0$, there exist $r_1, r_2>0$, such that $Z_1-r_1$ and $-Z_2-r_2$  are m-$\psi$-accretive operators with $\psi=\pi/4+\varepsilon$.
	
	In particular,  $-Z_1+r_1$ and $Z_2+r_2$   generates  holomorphic  $C_0$-semigroup of contraction operators $\mathcal{T}_1(z)$ and $\mathcal{T}_2(z)$ of angle $\dfrac{\pi}{2}-\psi$.

\end{proposition}
%%%%%%%%%%%%%%%%%%%%%%%%%%%%%%%%%%%%%%%%%%%%%%%%%%%%%%%%%%%%%%%%%%%%

\begin{proof}
	  Now assume  that $\mathcal{D}(\Lambda^{\frac{1}{2}}) \subset \mathcal{D} (B) $. It follows that
	\begin{equation}\label{AL1sur2}
	\left\|Bx\right\| \leq  a  \left\|x\right\|+ b\left\|\Lambda^{\frac{1}{2}}x\right\|
	\end{equation}
	for all $x\in \mathcal{D} (\Lambda^{\frac{1}{2}})$ and for some nonnegative constants $a$ and $b$.
	On the other hand, since $\Lambda^{\frac{1}{2}}$ is m-$(\pi/4)$-accretive  and both $B$ and $-B$ satisfy \eqref{AL1sur2}, by \cite[Theorem IX-2.4]{Kato},   we obtain the desired results.
\end{proof}
%%%%%%%%%%%%%%%%%%%%%%%%%%%%%%%%%%%%%%%%%%%%%%%%%%%%%%%%%%%%%%%%%%%%%%%%%%%%%%%%%%%
\begin{remark}\label{a0} If $a=0$ in \eqref{AL1sur2}, we have $r_1=r_2=0$ (cf. \cite[Theorem IX-2.4]{Kato}).
\end{remark}
%%%%%%%%%%%%%%%%%%%%%%%%%%%%%%%%%%%%%%%%%%%%%%%%%%%%%%%%%%%%%%%
\begin{proposition}\label{prop3.8}  If  $\mathcal{D}(\Lambda^{\frac{1}{2}}) \subset \mathcal{D} (B) $ and $B$ is accretive, then   $Z_1$ is m-$\pi/4$-accretive operators.	In particular, $-Z_1$   generates  holomorphic  $C_0$-semigroup $\mathcal{T}_1(z)$  of angle $\pi/4$.
	
	Further,  if $B$ is $\theta$-accretive with $0\leq \theta<\pi/2$, then $\mathcal{N}(Z_1)=\mathcal{N}(B)\cap\mathcal{N}(\Lambda^{\frac{1}{2}})$.
\end{proposition}
%%%%%%%%%%%%%%%%%%%%%%%%%%%%%%%%%%%%%%%%%%%%%%%%%%%%%%%%%%%%%%%%%%%%
\begin{proof}By \cite[Theorem 6.10]{Pazy83}, we have for an arbitrary $\rho>0$, 
	\begin{equation}\label{lalmda120}
	\left\|\Lambda^{\frac{1}{2}}x\right\|^2\leq \dfrac{1}{\pi^2} (\rho\left\|x \right\|^2+\dfrac{1}{\rho}\left\|\Lambda x\right\|^2),
	\end{equation}
	for all $ x\in \mathcal{D}  (B^2)$. 
	
Thus	by \eqref{lalmda120} and \eqref{AL1sur2}, we obtain
	\begin{equation*}%\label{BCB}
	\left\|Bx\right\|^2 \leq  2a (1+ \dfrac{\rho}{\pi^2}) \left\|x\right\|^2+\dfrac{2b }{\pi^2\rho_2}\left\|\Lambda x\right\|^2,
	\end{equation*}
	for all $ x\in \mathcal{D}  (B^2)$ and an arbitrary  $\rho>0$. Thus
	\begin{equation*}%\label{BCB}
	\left\|B (t+ \Lambda^{\frac{1}{2}})^{-1}x\right\|^2 \leq  2a(1+ \dfrac{\rho}{\pi^2}) \left\|(t+ \Lambda^{\frac{1}{2}})^{-1}x\right\|^2+\dfrac{2b}{\pi^2\rho}\left\|\Lambda (t+ \Lambda^{\frac{1}{2}})^{-1}x\right\|^2,
	\end{equation*}
	for all $x\in \mathcal{H}$. 
	
	Hence
	\begin{equation*}%\label{BCB}
	\left\|B (t+ \Lambda^{\frac{1}{2}})^{-1}\right\|^2 \leq  \dfrac{2a}{t^2}(1+ \dfrac{\rho}{\pi^2}) +\dfrac{2b}{\pi^2\rho}\left\|\Lambda (t+ \Lambda^{\frac{1}{2}})^{-1}\right\|^2.
	\end{equation*}
	Letting $t$ to $+\infty$, we assert that
	$$	M=\sup_{t>0} \left\| B(t+ \Lambda^{\frac{1}{2}})^{-1} \right\|<\dfrac{2b}{\pi^2\rho^2}. $$
	(cf. \cite[Proposition 2.12]{Yoshikawa72}). Since $\rho$ is arbitrary, we  can  choose it such that $\dfrac{2b}{\pi^2\rho^2}<1$. Thus $Z_1$  is  m-accretive. Since $B$ is accretive and $\Lambda^{\frac{1}{2}}$ is m-$(\pi/4)$-accretive, then $Z_1$  is  m-$(\pi/4)$-accretive. By \cite[Theorem IX-1.24] {Kato},	the factor $-Z_1$   generates  holomorphic  $C_0$-semigroup $\mathcal{T}_1(z)$  of angle $\dfrac{\pi}{4}$.
	
	The inclusion $\mathcal{N}(B)\cap\mathcal{N}(\Lambda^{\frac{1}{2}}) \subset \mathcal{N} (Z_1)$ is  obvious. Conversely, let $x\in \mathcal{D}  (Z_1)$, $x\neq 0$, such that $Z_1x=0$,  we have
	$$Re< Z_1x,x>=Re< Bx,x>+Re <\Lambda^{\frac{1}{2}} x,x>,$$
	then 
	$$Re< Bx,x> \leq Re< Z_1x,x> \quad \text{ and } \quad Re <\Lambda^{\frac{1}{2}} x,x> \leq Re< Z_1x,x>.$$
	Therefore, $Re< Z_1x,x>=0$ implies that  $Re< Bx,x>=0$ and  $Re< \Lambda^{\frac{1}{2}}x,x>=0$. 
	On the other hand, we have
	$$\left|Im <Bx,x>\right| \leq \tan(\theta) Re< Bx,x> $$ and $$\left|Im <\Lambda^{\frac{1}{2}}x,x>\right| \leq Re< \Lambda^{\frac{1}{2}}x,x>.$$
	Thus,
	$$ Im <Bx,x>=0 \quad \text{ and } \quad  Im <\Lambda^{\frac{1}{2}}x,x>=0,$$
	hence
	$$   <Bx,x>=0 \quad \text{ and } \quad   <\Lambda^{\frac{1}{2}}x,x>=0 .$$
	Since $B$ is m-$\theta$-accretive (see Remark \ref{remark3.2}) and $\Lambda^{\frac{1}{2}}$ is m-$(\pi/4)$-accretive, we conclude that
	$$ Bx= 0 \quad \text{ and } \quad  \Lambda^{\frac{1}{2}} x= 0. $$
	Consequently,  $ \mathcal{N} (Z_1)\subset\mathcal{N} (B)\cap \mathcal{N} ( \Lambda^{\frac{1}{2}}).$
\end{proof}
%%%%%%%%%%%%%%%%%%%%%%%%%%%%%%%%%%%%%%%%%%%%%%%%%%%%%%%%%%%%%%%%%%%%%%%%%%%%%%%%%%%%
\begin{remark}
	In Proposition \ref{thm:1D1}, if we assume only $B$ is m-accretive, $C + B^2$ need not be m-accretive, because $B^2$ fails to be accretive (with the same vertex as $B$) even in the case of an accretive matrix $B$ with numerical
	range contained in a sector of angle less than $\pi/4$,  as the following example shows.
\end{remark}
%%%%%%%%%%%%%%%%%%%%%%%%%%%%%%%%%%%%%%%%%%%%%%%%%%%%%%%%%%%%%%%%%%
\begin{example}\label{b2b} Let $\mathcal{H}=\mathbb{C}^2$ and
	$$B=\begin{bmatrix}
	4-i& 4i\\
	4i & 16+4i\\
	\end{bmatrix}.$$
	
For $x=(x_1 , x_2)\in \mathbb{C}^2$; we have
$$ Re<Bx,x>=4\left|x_1 \right|^2+16 \left|x_2 \right|^2$$
and
\begin{align*}
Im<Bx,x>& =-\left|x_1 \right|^2+8 Re(x_1 \overline{x_2})+ 4 \left|x_2 \right|^2\\
&\leq 3\left|x_1 \right|^2+8 \left|x_2 \right|^2\\
& < Re<Bx,x> .
\end{align*}
Thus
$$W(B)\subsetneq S_{\pi/4}.$$
However, for $x=(1,0)$, we have
 $$<B^2x, x>=-1-8i,$$ 
it follows that $W(B^2)$ is not a subset of the right half complex plane.
\end{example}
%%%%%%%%%%%%%%%%%%%%%%%%%%%%%%%%%%%%%%%%%%%%%%%%%%%%%%
\begin{remark}The  operator pencil $\mathcal{Q}$ is not necessarily an accretive, because we can find an  eigenvalues  not  located in the closed right  half-plane. Indeed, let $\lambda$ be an eigenvalue of  $\mathcal{Q}$ and  $v\in \mathcal{D} (\mathcal{Q})$  its corresponding  eigenvector with $\left\| v\right\| =1$. Let us remark that if $\lambda =0$, then $0=<Cv,v>$ and hence $0\in W(C)$. In the sequel we assume that $\lambda \neq 0$ with $\lambda=\alpha+ i \beta $. Then 
		$$<\mathcal{Q}(\lambda)v,v>=0, $$
		and consequently, taking real and imaginary parts,
		$$(\alpha^2-\beta^2)-2\alpha Re<Bv,v>+2\beta Im<Bv,v>-Re<Cv,v>=0, $$
		and
		$$
		2\alpha\beta-2\beta Re<Bv,v>-2\alpha Im<Bv,v>-Im<Cv,v>=0.
		$$ 
		It follows that
		$$Re<Cv,v>=(\alpha^2-\beta^2)-2\alpha Re<Bv,v>+2\beta Im<Bv,v>,$$
		and
		$$
		Im<Cv,v>=2\alpha\beta-2\beta Re<Bv,v>-2\alpha Im<Bv,v>.
		$$ 
		Since $Re<Cv,v> \geq 0$, we obtain from the first relation,
		$$2\alpha Re<Bv,v>\leq \alpha^2-\beta^2+2\beta Im<Bv,v>.$$
		The fact that $\left| Im<Bv,v>\right| \leq Re<Bv,v>$, we get
		$$2\alpha Re<Bv,v>\leq \alpha^2-\beta^2+2 \left| \beta\right| Re<Bv,v>.$$
		Thus
		$$2(\alpha - \left| \beta\right|) Re<Bv,v>\leq \alpha^2-\beta^2. $$
		Now, if assume  $\left| \alpha \right| \leq \left| \beta \right|$, it follows that
		$$(\alpha - \left| \beta\right|) Re<Bv,v>\leq 0. $$
		Consequently, $ \alpha  \leq \left| \beta \right|$.
\end{remark}
%%%%%%%%%%%%%%%%%%%%%%%%%%%%%%%%%%%%%%%%%%%%%%%%%%%%
\begin{remark}
 Similar results  can be obtained by interchanging the role of $B^2$ and $C$, be careful with domains. In this case we have $\mathcal{D} (C)\subset  \mathcal{D} (B^2)\subset \mathcal{D} (B)$. 
\end{remark}
%%%%%%%%%%%%%%%%%%%%%%%%%%%%%%%%%%%%%%%%%%%%%%%%%%%%
%%%%%%%%%%%%%%%%%%%%%%%%%%%%%%%%%%%%%%%%%%%%%%%%%%%%%
\section{An application to an abstract second order differential equation}
%%%%%%%%%%%%%%%%%%%%%%%%%%%%%%%%%%%%%%%%%%%%%%%%%%%%%%%%%%%%%
Let us consider, in the complex Hilbert space $\mathcal{H}$,  the following  abstract second order differential equation
\begin{equation}\label{2EDA}
u^{\prime \prime }(x)-2Bu^{\prime }(x)-Cu(x)=f(x),  \qquad x \in(0, 1),
\end{equation}
under the boundary conditions
\begin{equation}\label{2EDABC}
  u(0)=u_0, \qquad  u(1)=u_1,
\end{equation}
where  $B$ and $C$  are  two closed  operators in a Hilbert space with domains $\mathcal{D} (B)$ and $\mathcal{D} (C)$, respectively, $f\in L^p(0,1; \mathcal{H})$, $1<p<\infty$ and $u_0$, $u_1$ are given elements in $\mathcal{H}$. We  seek for a strict solution $u$ to \eqref{2EDA}-\eqref{2EDABC}, i.e. a function $u$ such that
\begin{equation}
\left \{
\begin{array}{l}
i)\text{ }u\in W^{2,p}(0,1; \mathcal{H})\cap L^{p}(0,1; \mathcal{D} (C))),\text{ } u^{\prime }\in
L^{p}(0,1; \mathcal{D} (B)), \\
ii)\text{ }u \text{ satisfies } \eqref{2EDA}-\eqref{2EDABC}.%
\end{array}%
\right.  \label{OptimalRegularities}
\end{equation}

%%%%%%%%%%%%%%%%%%%%%%%%%%%%%%%%%%%%%%%%%%%%%%%%%%%%%%%%%%
\begin{theorem}\label{thm:solEDA} Let $B$   and  $C$  two operators in a Hilbert space $\mathcal{H}$ such that
	\begin{enumerate}
		\item $B^2$  is m-accretive and  $C$  is accretive     satisfy one of conditions of Proposition \ref{thm:1D1}.
		\item  $ \mathcal{D} ((B^2 +C)^{\frac{1}{2}})\subset \mathcal{D}  (B)$.
		\item   $B(\mathcal{D}  (B^2))\subset \mathcal{D}  (B^2) $ and $(B^2 +C)^{\frac{1}{2}}(\mathcal{D}  (B^2))\subset \mathcal{D}  (B^2) $.
		\item $(B^2 +C)^{-\frac{1}{2}}$ exist and bounded.
		\item $B(B^2 +C)^{\frac{1}{2}}=(B^2 +C)^{\frac{1}{2}}B$ on $\mathcal{D} (B^2)$.
		\item  $f\in L^{p}(0,1; \mathcal{H})$ with $1<p<\infty $.
	\end{enumerate}
 Then the problem \eqref{2EDA}-\eqref{2EDABC} has a classical solution $u$ if and only if 
\begin{equation*}
Z_1^2e^{. -Z_1}u_0, \quad Z_1^2e^{. - Z_1}u_1\in L^{p}(0,1; \mathcal{H}) .
\end{equation*}%
In this case, $u$ is uniquely determined by 
\begin{align*}
u(x)& =(I-e^{Z_2-Z_1})^{-1}e^{xZ_2}u_0 +(I-e^{Z_2-Z_1})^{-1}e^{-(1-x)Z_1}u_1\\
     & - (I-e^{Z_2-Z_1})^{-1}e^{xZ_2}e^{-Z_1}\left(u_1 -(Z_2-Z_1)^{-1} \int_{0}^{1}e^{\left(1-s\right)
     	Z_2}f(s)ds\right) \\
     & - (I-e^{Z_2-Z_1})^{-1}e^{-(1-x)Z_1}e^{Z_2}\left(u_0 -(Z_2-Z_1)^{-1} \int_{0}^{1}e^{-sZ_1}f(s)ds\right) \\
     &-(I-e^{Z_2-Z_1})^{-1}(Z_2-Z_1)^{-1}e^{xZ_2}\int_{0}^{1}e^{-sZ_1}f(s)ds\\
     &+ (I-e^{Z_2-Z_1})^{-1}(Z_2-Z_1)^{-1}e^{-(1-x)Z_1} \int_{0}^{1}e^{-\left(1-s\right)
     	Z_2}f(s)ds\\
     & +(Z_2-Z_1)^{-1}\int_{0}^{x}e^{\left(	x-s\right) Z_2}f(s)ds+(Z_2-Z_1)^{-1}\int_{x}^{1}e^{\left(	x-s\right) Z_1}f(s)ds.
\end{align*}
\end{theorem}
%%%%%%%%%%%%%%%%%%%%%%%%%%%%%%%%%%%%%%%%%%%%%%%%%%%%%%%%%%%%%
\begin{proof}  Under the assumptions, by Proposition \ref{thm:1D1S} and Remark \ref{a0}, the factors $-Z_1$ and $Z_2$   generates bounded  holomorphic  $C_0$-semigroup $(e^{-tZ_1})_{t\geq 0}$ and $(e^{tZ_2})_{t\geq 0}$, respectively. Also, $\mathcal{D} (Z_1)=\mathcal{D} (Z_2)=\mathcal{D} (\Lambda^{1/2})$ and
	$$ \mathcal{D} (Z_1Z_2)=\left\lbrace x \in\mathcal{D} (Z_2) ; Z_2x\in \mathcal{D} (Z_1)  \right\rbrace = \left\lbrace x \in\mathcal{D} (Z_2) ; Z_2x\in \mathcal{D} (Z_2)  \right\rbrace= \mathcal{D} (Z_2^2),$$ 
	$$ \mathcal{D} (Z_2Z_1)=\left\lbrace x \in\mathcal{D} (Z_1) ; Z_1x\in \mathcal{D} (Z_2)  \right\rbrace = \left\lbrace x \in\mathcal{D} (Z_1) ; Z_1x\in \mathcal{D} (Z_1)  \right\rbrace= \mathcal{D} (Z_1^2). $$ 
	But
	$$\mathcal{D} (Z_1^2)=\left\lbrace x \in\mathcal{D} (\Lambda^{1/2}) ; Z_1x\in \mathcal{D} (\Lambda^{1/2})  \right\rbrace $$
	and
	$$\mathcal{D} (Z_2^2)=\left\lbrace x \in\mathcal{D} (\Lambda^{1/2}) ; Z_2x\in \mathcal{D} (\Lambda^{1/2})  \right\rbrace.$$
	The fact that, $B(\mathcal{D}  (B^2))\subset \mathcal{D}  (B^2) $ and $\Lambda^{\frac{1}{2}}(\mathcal{D}  (B^2))\subset \mathcal{D}  (B^2) $, we obtain $\mathcal{D}  (B^2)\subset \mathcal{D} (Z_1^2)$ and $\mathcal{D}  (B^2)\subset \mathcal{D} (Z_2^2)$, with $\mathcal{D}  (B^2)$ densely defined on $\mathcal{H}$. Furthermore, $e^{-tZ_1}u_0\in \mathcal{D} (Z_1^n)$ and $e^{tZ_2}u_1\in \mathcal{D} (Z_2^n)$ for all $u_0$ ,$u_1\in\mathcal{H}$, $t > 0$ and $n\in\mathbb{N}$. Hence $u(x)\in \mathcal{D} (C)$ for all $x\in (0, 1)$. Since the two $C_0$-semigroups are holomorphic, $u(.)$ can be differentiated any numbers of times.  Now, by taking $-B$ instead $B$, $A=-C$, $L=-Z_1$ and $M=Z_2$ in \cite[Theorem 5.]{Favini2008}, all assumptions of this theorem are fulfilled. Hence we obtain the desired result.
\end{proof}

%%%%%%%%%%%%%%%%%%%%%%%%%%%%%%%%%%%%%%%%%%%%%%%%%%%%%%%%%%%%%
\section{An example of a second order partial differential equation}%
%%%%%%%%%%%%%%%%%%%%%%%%%%%%%%%%%%%%%%%%%%%%%%%%%%%%%%%%%%%%%
%%%%%%%%%%%%%%%%%%%%%%%%%%%%%%%%%%%%%%%%%%%%%%%%%%%%%
The aim of  this section is to use the obtained results to discuss the existence, uniqueness, and maximal regularity of the strict solution  for the following non-homogeneous second order differential equation,
%%%%%%%%%%%%%%%%%%%%%%%%%%%%%%%%%%%%%%%%%%%%%%%%%%%%%%%%%%%%%%%%%%%%%%%%%%%%%%%%%%%%%%%%%%%%%%%%%%%%%%%%%%%
\begin{equation*}
(E)\left \{
\begin{array}{l}
\dfrac{\partial ^{2}u}{\partial x^{2}}(x,y)-2  p_0(y)\dfrac{\partial ^{2}u}{%
	\partial y\partial x}(x,y)-2  p_1(y)\dfrac{\partial u}{
	\partial x}(x,y)+\alpha  p_0(y)\dfrac{\partial u}{
	\partial y}(x,y)\\
\qquad \qquad \qquad \qquad +(\alpha p_1(y)+\beta  )u(x,y)-\gamma u(x,y)=f(x,y),\quad x\in (0,1),\quad y \in (0,1)
\\
\\
u(0,y)=u_0(y),\quad u(1,y)=u_1(y),\quad y\in (0,1)\\
u(x,0)=u(x,1)=0,\quad x\in (0,1)\\
\dfrac{\partial u}{\partial x}(x,0)=\dfrac{\partial u}{\partial x}(x,1)=0,\quad x\in (0,1)
\end{array}%
\right.
\end{equation*}
where, 
\begin{itemize}
	\item $f\in L^p(0,1; L^2(0,1; \mathbb{C}))$, $1<p<\infty$, 
	\item  $\alpha\in \mathbb{R}$, $\beta\in \mathbb{C}$, $p_0 , p_1 \in C^1(0,1)$ and $p_0(x)\neq 0 $ for all $x\in [0,1]$.
	\item  $\gamma =-(\dfrac{r+1}{4\varepsilon}M_1+M_2)$, with  $r>0$ and $\varepsilon$ are arbitrary and  chosen such that $m_0-\varepsilon(1+ r)M_1>0$, for some nonegative constants $m_0$, $M_1$ and $M_2$ are described below.
\end{itemize}
%%%%%%%%%%%%%%%%%%%%%%%%%%%%%%%%%%%%%%%%%%%%%%%%%%
The second order differential equation $(E)$
is equivalent to
%%%%%%%%%%%%%%%%%%%%%%%%%%%%%%%%%%%%%%%%%%%%%%%%%%%%%
\begin{equation}\label{2EDAuxy}
\dfrac{\partial ^{2}u}{\partial x^{2}}(x,y)-2  B\dfrac{\partial u}{
	\partial x}(x,y)+Cu(x,y)-\gamma u(x,y)
=f(x,y),\quad x\in (0,1),\quad y \in (0,1).
\end{equation}
%%%%%%%%%%%%%%%%%%%%%%%%%%%%%%%%%%%%%%%%%%%%%%%%%%%%%
with the boundary conditions
\begin{equation}\label{2EDAuxyBC}
u(0,y)=u_0(y),\quad u(1,y)=u_1(y),\quad y\in (0,1),
\end{equation}
where,
%%%%%%%%%%%%%%%%%%%%%%%%%%%%%%%%%%%%%%%%%%%%%%%%%%%%%
\begin{equation*}
\left \{
\begin{array}{l}
B=  p_0\dfrac{\partial }{%
	\partial y} +  p_1,\quad \mathcal{D} (B)=\{ \psi  \in H^1(0,1) : \psi(0)=\psi(1)=0 \}
\\
\text{ and }\\
C= \alpha p_0\dfrac{\partial }{%
	\partial y}+ (\alpha  p_1 +\beta),\quad \mathcal{D} (C)=\{ \phi \in H^1(0,1) : \phi(0)=\phi(1)=0 \}.
\end{array}%
\right.
\end{equation*}
%%%%%%%%%%%%%%%%%%%%%%%%%%%%%%%%%%%%%%%%%%%%%%%%%%%%
with $ \phi(y)= u(x,y)$ and $\psi(y)=\dfrac{\partial u}{
	\partial x}(x,y)$, $x\in (0,1),\quad y \in (0,1)$.
%%%%%%%%%%%%%%%%%%%%%%%%%%%%%%%%%%%%%%%%%%%%%%%%%%%%%%%%%%%%%%%%
%%%%%%%%%%%%%%%%%%%%%%%%%%%%%%%%%%%%%%%%%%%%%%%%%%%%%
We  seek for a strict solution $u(. , y)$ to \eqref{2EDAuxy}-\eqref{2EDAuxyBC}, i.e. a function $u(. , y)\in  L^2(0,1; \mathbb{C})$ such that \ref{OptimalRegularities} holds. This will be done  by the following preparatory results. 
%%%%%%%%%%%%%%%%%%%%%%%%%%%%%%%%%%%%%%%%%%%%%%%%%%%%%%%%%%

\textbf{Claim 1.} \textit{The operator $-B^2-\gamma I$ is  m-$\omega$-accretive, with $\omega=\arctan(\dfrac{1}{r})$.}

%%%%%%%%%%%%%%%%%%%%%%%%%%%%%%%%%%%%%%%%%%%%%%%%%%%%%%%%%%%%
\begin{proof}   For $\psi \in \mathcal{D} (B^2)\subset \{ \psi  \in H^2(0,1) : \psi(0)=\psi(1)=0 \}\subset \mathcal{D} (B)$, we have
	$$-B^2\psi=\varphi_0 \psi''+\varphi_1 \psi'+\varphi_2\psi,$$
	with $\varphi_0 =-p_0^2$, $\varphi_1 =-p_0(p'_0+2p_1)$ and $\varphi_2 =-(p_1^2+p_0p'_1).$ 
	Under the assumptions there exists a nonegative constants $m_0$, $M_0$ and $M_1$ such that
\begin{equation}\label{m0M1M2}
 -\varphi_0>m_0>0, \qquad \qquad  \left| \varphi_1-\varphi'_0\right| \leq  M_1,  \qquad \text{ and }\qquad  \left|  \varphi_2\right|\leq M_2.
\end{equation}
	By  \cite[Example V-3.34]{Kato}, $-B^2$ is m-$\omega$-accretive operator with vertex $\gamma$, where $\gamma =-(\dfrac{r+1}{4\varepsilon}M_1+M_2)$, $\omega=\arctan(\dfrac{1}{r})$, $r>0$ and $\varepsilon$ is chosen such that $m_0-\varepsilon(1+r)M_1>0$. Hence the operator $-B^2-\gamma I$ is  m-$\omega$-accretive, with $\omega=\arctan(\dfrac{1}{r})$.%Further,  $\Lambda$ generates  holomorphic  $C_0$-semigroup of contraction operators $\mathcal{T}(z)$ of angle $\dfrac{\pi}{2}-\omega$, see \cite[Example IX-1.26]{Kato}. It is not difficult to prove $BC=CB$ on  $\mathcal{D} (B)$.
\end{proof}
%%%%%%%%%%%%%%%%%%%%%%%%%%%%%%%%%%%%%%%%%%%%%%%%%%%%%%%%%%%

\textbf{Claim 2.} \textit{ If $\alpha  p_1+Re(\beta)- \dfrac{\alpha}{2}p'_0 \geq 0$ then $C$ is an  accretive operator.}
  
%%%%%%%%%%%%%%%%%%%%%%%%%%%%%%%%%%%%%%%%%%%%%%%%%%%%%%%%%%%%
%%%%%%%%%%%%%%%%%%%%%%%%%%%%%%%%%%%%%%%%%%%%%%%%%%%%%%%%%%%%%%%%%%%
\begin{proof}Let $ \psi\in \mathcal{D}  (C)$, we have
	$$<C\psi, \psi>=\alpha \int_0^1 p_0(y)  <\psi'(y),\overline{\psi(y)}>dy+ \int_0^1 (\alpha  p_1(y)+\beta)  \left| \psi(y)\right|^2 dy.$$
	By integration by parts,
	$$<C\psi, \psi>=-\alpha \int_0^1 p_0(y) <\psi(y),\overline{\psi'(y)}>dy+ \int_0^1 ( \alpha  p_1(y)+\beta -\alpha p'_0(y)) \left| \psi(y)\right|^2 dy.$$
	Also
	$$\overline{<C\psi, \psi>}=<\psi, C\psi>=\alpha \int_0^1 p_0(y) <\psi(y),\overline{\psi'(y)}>dy+ \int_0^1 (\alpha  p_1(y)+\overline{\beta})  \left| \psi(y)\right|^2 dy.$$
	Thus
	$$Re<C\psi, \psi>= \int_0^1 (\alpha  p_1+Re(\beta)- \dfrac{\alpha}{2}p'_0) \left| \psi(y)\right|^2 dy.$$
	Hence the desired result.
\end{proof}
%%%%%%%%%%%%%%%%%%%%%%%%%%%%%%%%%%%%%%%%%%%%%%%%%%%%%%%%%%%%%%%%%%
If we take $\alpha =1$ and $\beta =0$ in Claim 2, we obtain

\textbf{Claim 3.} \textit{ If $p_1- \dfrac{1}{2}p'_0\geq 0$ then $B$ is an  accretive operator. In particular, by Remark \ref{remark3.2}, $B$ is m-accretive.}

%%%%%%%%%%%%%%%%%%%%%%%%%%%%%%%%%%%%%%%%%%%%%%%%%%%%%%%%%%%

\textbf{Claim 4.} \textit{  If $p_1- \dfrac{1}{2}p'_0\geq 0$  and $\alpha  p_1+Re(\beta)- \dfrac{\alpha}{2}p'_0 \geq 0$, then $-\Lambda=-B^2 +C -\gamma I$ with domain $ \mathcal{D}  (B^2)$  is m-accretive. Also,  $-\Lambda$ admits an unique square root  $(-\Lambda)^{1/2}$  m-$(\pi/4)$-accretive.}

%%%%%%%%%%%%%%%%%%%%%%%%%%%%%%%%%%%%%%%%%%%%%%%%%%%%%%%%%%%%

%%%%%%%%%%%%%%%%%%%%%%%%%%%%%%%%%%%%%%%%%%%%%%%%%%%%%%%%%%%%%%%%%%%
\begin{proof} By Claim 1. $-B^2 -\gamma I$ with domain $ \mathcal{D}  (B^2)$  is m-accretive, by Claim 2.  $C$ is an  accretive and  by Claim 3. $B$ is an  accretive operator. Also, $ \mathcal{D}  (B)=\mathcal{D}  (C)$.  Now the desired result holds from the \textbf{(C.4)} and Proposition \ref{thm:1D1}.
\end{proof}
%%%%%%%%%%%%%%%%%%%%%%%%%%%%%%%%%%%%%%%%%%%%%%%%%%%%%%%%%%%%%%%%%%

%%%%%%%%%%%%%%%%%%%%%%%%%%%%%%%%%%%%%%%%%%%%%%%%%%%%%%%%%%%

\textbf{Claim 5.} \textit{ If $p''_0$ is continuous on $[0 ,1]$, then $(-B^2+C-\gamma I)^{-1}$  exist and bounded.}

%%%%%%%%%%%%%%%%%%%%%%%%%%%%%%%%%%%%%%%%%%%%%%%%%%%%%%%%%%%%
%%%%%%%%%%%%%%%%%%%%%%%%%%%%%%%%%%%%%%%%%%%%%%%%%%%%%%%%%%%%
\begin{proof}   As before; for $\psi \in \mathcal{D} (B^2)\subset \{ \psi  \in H^2(0,1) : \psi(0)=\psi(1)=0 \}\subset \mathcal{D} (B)$, we have
		$$[-B^2+C-\gamma I]\psi=\varphi_0 \psi''+(\varphi_1+\alpha p_1) \psi'+(\varphi_2+\alpha p_1+\beta-\gamma) \psi,$$
	with $\varphi_0 =-p_0^2$, $\varphi_1 =-p_0(p'_0+2p_1)$ and $\varphi_2 =-(p_1^2+p_0p'_1).$ Since $p''_0$ and $p'_1$ are continuous on $[0 ,1]$, it follows that $\varphi''_0$, $\varphi'_1+\alpha p'_1$ and $\varphi_2+\alpha p_1+\beta-\gamma$ are are continuous on $[0 ,1]$. By a similar way as in  \cite[Section 3-III. p. 146-149]{Kato}, we prove that $(-B^2+C -\gamma I)^{-1}$  exist and bounded.
\end{proof}
%%%%%%%%%%%%%%%%%%%%%%%%%%%%%%%%%%%%%%%%%%%%%%%%%%%%%%%%%%%
 Combining Claim 4.  Corollary \ref{root}, Proposition \ref{thm:1D1S} and Proposition \ref{prop3.8}, we obtain,

\textbf{Claim 6.} \textit{	The operators $Z_1=B+(-\Lambda)^{1/2}$ and $Z_2=B-(-\Lambda)^{1/2}$ with domain $\mathcal{D}(\Lambda^{\frac{1}{2}})\subset \mathcal{D} (B)$ are $B^2$-bounded  and closed operators. Furthermore, the closure of the restriction of  $Z_i$ to  $\mathcal{D}  (B^2)$ is again $  Z_i$, $i=1,2$, $-Z_1$   generates  holomorphic  $C_0$-semigroup $\mathcal{T}_1(z)$  of angle $\pi/4$ and $Z_1+r$   generates  holomorphic  $C_0$-semigroup $\mathcal{T}_1(z)$  of angle $\pi/4+\varepsilon$, for some $\varepsilon>0$ and $r>0$. }

We are now ready to state the following existence and uniqueness  result.

%%%%%%%%%%%%%%%%%%%%%%%%%%%%%%%%%%%%%%%%%%%%%%%%%%%%%%%%%%%%%%%%%%%%%%
%%%%%%%%%%%%%%%%%%%%%%%%%%%%%%%%%%%%%%%%%%%%%%%%%%%%%%%%%%
\begin{theorem}Let the equation $(E)$ on  $\mathcal{H}=L^2(0,1; \mathbb{C}) $. Assume that
	\begin{enumerate}
		\item $f\in L^p(0,1; \mathcal{H})$, $1<p<\infty$, 
		\item  $\alpha\in \mathbb{R}$, $\beta\in \mathbb{C}$, $p_0 \in C^2(0,1)$ , $p_1 \in C^1(0,1)$ and $p_0(x)\neq 0 $ for all $x\in [0,1]$,
		\item $p_1- \dfrac{1}{2}p'_0\geq 0$ and $\alpha ( p_1- \dfrac{1}{2}p'_0)+Re(\beta) \geq 0$,
		\item  $\gamma =-(\dfrac{r+1}{4\varepsilon}M_1+M_2)$, with  $r>0$ and $\varepsilon$ are arbitrary and  chosen such that $m_0-\varepsilon(1+r)M_1>0$, for some nonegative constants $m_0$, $M_1$ and $M_2$ are  given by \eqref{m0M1M2}.
		\item $B(-\Lambda)^{1/2}=(-\Lambda)^{1/2}B$ on $\mathcal{D} (B^2)$.
	\end{enumerate}
	Then the problem \eqref{2EDAuxy}-\eqref{2EDAuxyBC} has a classical solution $u$ if and only if 
	\begin{equation*}
	Z_1^2e^{. -Z_1}u_0, \quad Z_1^2e^{. -Z_1}u_1\in L^{p}(0,1; \mathcal{H}) .
	\end{equation*}%
	In this case, $u$ is uniquely determined as in Theorem \ref{thm:solEDA}.
\end{theorem}
%%%%%%%%%%%%%%%%%%%%%%%%%%%%%%%%%%%%%%%%%%%%%%%%%%%%%%%%%%%%%
\begin{proof}
Thus the restriction of  $Z_1$  and $-Z_2$  to  $\mathcal{D}  (B^2)$  are m-$(\pi/4)$-accretive  operators. Also; by Claim 4., the inverse of $(-\Lambda)^{1/2}$ exist and bounded. Thus, all assumptions of Theorem \ref{thm:solEDA} are fulfilled. Consequently, we get the desired result.

\end{proof}

%%%%%%%%%%%%%%%%%%%%%%%%%%%%%%%%%%%%%%%%%%%%%%%%%%%%%%%%.

%%%%%%%%%%%%%%%%%%%%%%%%%%%%%%%%%%%%%%%%%%%%%%%%%%%%%%%%

%%%%%%%%%%%%%%%%%%%%%%%%%%%%%%%%%%%%%%%%%%%%%%%%%

\begin{thebibliography}{00}
%%%%%%%%%%%%%%%%%%%%%%%%%%%%%%%%%%%%%%%%%%%%%%%%%%%%%%%%%%%%%
%%%%%%%%%%%%%%%%%%%%%%%%%%%%%%%%%%%%%%%%%%%%%%%%%%%%%%%%%%%%%%%%
\bibitem{AglazinKiiko1}  S. D. Aglazin, I. A. Kiiko, \textit{Numerical--analytic investigation of the flutter of a panel
	of arbitrary shape in a desighn}. J. Appl. Math. Mech. \textbf{61} (1997), 171--174
%%%%%%%%%%%%%%%%%%%%%%%%%%%%%%%%%%%%%%%%%%%%%%%%%%%%%%%%%%%%%%%%%%%%%%%%%
\bibitem{Arlinskii2010} Y. M.  Arlinski\u{\i}, V. Zagrebnov, \textit{Numerical range and quasi-sectorial contractions}.  J. Math. Anal. Appl. 366 (2010); 33-43.
%%%%%%%%%%%%%%%%%%%%%%%%%%%%%%%%%%%
\bibitem{Artamonov1} N. Artamonov, \textit{Estimates of solutions of certain classes of second-order differential
	equaitons in a Hilbert space}. Sbornik Mathematics, \textbf{194}:8 (2003), 1113--1123
%%%%%%%%%%%%%%%%%%%%%%%%%%%%%%%%%%%%%%%%%%%%%
\bibitem{Balakrishnan60} A. V. Balakrishnan, \textit{Fractional powers of closed operators and the semi-groups generated by them}, Pacific J. Math. 10 (1960), 419--437.
%%%%%%%%%%%%%%%%%%%%%%%%%%%%%%%%%%%%%%%%%%%%%%%
%%%%%%%%%%%%%%%%%%%%%%%%%%%%%%%%%%%%%%%%%%%%%%%
\bibitem{israel} A. Ben-Israel and T. N. E. Greville, \emph{Generalized Inverses: Theory and Applications}.
Second edition, Springer, New York, 2003.
%%%%%%%%%%%%%%%%%%%%%%%%%%%%%%%%%%%%%%%%%%%%%%%%%%
%%%%%%%%%%%%%%%%%%%%%%%%%%%%%%%%%%%%%%%%%%%%%%%%%%%%%%%
%\bibitem{Chernoff72} P. R. Chernoff, \textit{Perturbations of dissipative operators with relative bound one}, Proc.
%Amer. Math. Soc. 33 (1972), 72--74.
%%%%%%%%%%%%%%%%%%%%%%%%%%%%%%%%%%%%%%%%%%%%%%%%%%
\bibitem{Duffin55} R. J. Duffin, \textit{A minimax theory for overdamped networks}, J. Rational Mech.
Anal. 4 (1955), 221--233.
%%%%%%%%%%%%%%%%%%%%%%%%%%%%%%%%%%%%%%%%%%%%%%%%%%%%%%
%\bibitem{Engel95} K-J. Engel, \textit{On perturbations of linear m-accretive operators on reflexive Banach spaces}, Mh. Math. 119 (1995), 259--265. 
%%%%%%%%%%%%%%%%%%%%%%%%%%%%%%%%%%%%%%%%%%%%%
\bibitem{Langer2004} D. Eschw\'e  M. Langer,\textit{ Variational principles for eigenvalues of selfadjoint
operator functions}, Integral Equations Operator Theory 49 (2004), 287--321.
%%%%%%%%%%%%%%%%%%%%%%%%%%%%%%%%%%%%%%%%%%%%%%%%%%%%%%%
\bibitem{Favini2008} A. Favini, R. Labbas, S. Maingot, H. Tanabe, A. Yagi,  \textit{A Simplified Approach in the Study of Elliptic Differential Equations in UMD Spaces and New Applications}.  Funkcialaj Ekvacioj,  vol. 51 (2008), pp. 165--187.
%%%%%%%%%%%%%%%%%%%%%%%%%%%%%%%%%%%%%%%%%%%%%%%%%%
\bibitem{Gustafson66} K. Gustafson,\textit{ A perturbation lemma}, Bull. Am. Math. Soc., 72 (1966), 334--338.
%%%%%%%%%%%%%%%%%%%%%%%%%%%%%%%%%%%%%%%%%%%%%%%%%%%%%%%%
%%%%%%%%%%%%%%%%%%%%%%%%%%%%%%%%%%%%%%%%%%%%%%%%%%%%%%%%
\bibitem{GustafsonR97} K. Gustafson, D. Rao, \emph{Numerical Range, the Field of Values of Linear Operators
	and Matrices}, Springer, New York, 1997.
%%%%%%%%%%%%%%%%%%%%%%%%%%%%%%%%%%%%%%%%%%%%%%%%%%%%%%%%%%
%\bibitem{GustafsonR77} K. Gustafson, D. Rao, \emph{Numerical range and accretivity of operator products},
%J. Math. Anal. Appl.  60 (3),  1977,  693-702.
%%%%%%%%%%%%%%%%%%%%%%%%%%%%%%%%%%%%%%%%%%%%%%%%%%%%%%%%%%
\bibitem{Hayashi2017} M. Hayashi and T. Ozawa, \textit{On Landau-Kolmogorov inequalities for dissipative operators}, Proc. Amer. Math. Soc. 145 (2017), 847-852.
%\bibitem{HessK70}  P. Hess, T. Kato, \textit{Perturbation of closed operators and their adjoints}, Comment. Math. Helv.
%45 (1970) 524--529.
%%%%%%%%%%%%%%%%%%%%%%%%%%%%%%%%%%%%%%
\bibitem{IlyushinKiiko1} A. A. Ilyushin, I. A. Kiiko, \textit{Vibrations of a rectangle plate in a supersonic aerodynamics and
	the problem of panel flutter}. Moscow Univ. Mech, Bull., \textbf{49} (1994), 40--44
%%%%%%%%%%%%%%%%%%%%%%%%%%%%%%%%%%%%%%%%%%%%%%%%%%%%%%%%%%%%%%%
\bibitem{KreinL64} M.G. Krein and H. Langer, \textit{On the theory of quadratic pencils of self-adjoint operators},
Dokl. Akad. Nauk SSSR 154 (1964), 1258--1261 (Russian); English transl.,
Soviet Math. Dokl. 5 (1964), 266--269.
%%%%%%%%%%%%%%%%%%%%%%%%%%%%%%%%%%%%%%%%%%%%
%\bibitem{Krol2009} S. Krol, \textit{Perturbation theorems for holomorphic semigroups}, J. Evol. Equ. 9 (2009), 449--468.
%%%%%%%%%%%%%%%%%%%%%%%%%%%%%%%%%%%%%%%%%%%%%%%%%%%%%%%%%%%%%%%%%%
%%%%%%%%%%%%%%%%%%%%%%%%%%%%%%%%%%%%%%%%%%%%%
\bibitem{Kato} T. Kato,
{\emph{Perturbation theory for linear operators}},
   Springer-Verlag, New York (1995).
%%%%%%%%%%%%%%%%%%%%%%%%%%%%%%%%%%%%%%%%%%%%%%%%%%%%%%%%%%%%%
%\bibitem{Kato60} T. Kato, \textit{Note on fractional powers of linear operators},  Proc. Japan Acad. 36 (1960), no. 3, 94--96. 
%%%%%%%%%%%%%%%%%%%%%%%%%%%%%%%%%%%%%%%%%%%%%%%%%%%%%%%%%%%%%%%
\bibitem{Kato61} T. Kato, \textit{Fractional powers  of dissipative operators},  Proc. Japan Acad. 13 (3) (1961), 246--274. 
%%%%%%%%%%%%%%%%%%%%%%%%%%%%%%%%%%%%%%%%%%%%%%%%%%%%%%%%%%%%%%%
\bibitem{Kato71} T. Kato, \textit{On an Inequality of Hardy, Littlewood, and Po1ya},  Advances in Math. 7  (1971), 217--218. 
%%%%%%%%%%%%%%%%%%%%%%%%%%%%%%%%%%%%%%%%%%%%%%%%%%%%%%%%%%%%%%%
\bibitem{Markus88}
A.~S. Markus, {\em Introduction to the spectral theory of polynomial operator
  pencils}, volume~71 of {\em Translations of Mathematical Monographs}. American Mathematical Society, Providence, RI, 1988.
%%%%%%%%%%%%%%%%%%%%%%%%%%%%%%%%%%%%%%%%%%%%%%%%%%%%%%%%%%%%%%%%%%%%%%%%%%%

%\bibitem{Melinkova83} I. V. Melnikova, \textit{On the Cauchy problem for a second-order equation}, Differ. equations,  vol. 19 no. 3 (1983), 537--538. [In Russian].
%%%%%%%%%%%%%%%%%%%%%%%%%%%%%%%%%%%%%%%%%%%%%%%%%%%%%%%%%%%%%
\bibitem{Moller2015} M. M\"{o}ller,  V. Pivovarchik, \emph{Spectral Theory of Operator Pencils, Hermite-Biehler Functions, and their Applications}, Birkh\"{a}user (2015).
%%%%%%%%%%%%%%%%%%%%%%%%%%%%%%%%%%%%%%%%%%%%%%%%%%%%%%%%%%%%%%%%%%%%%%%
\bibitem{Pazy83} A. Pazy,  \emph{Semigroups of Linear Operators and Applications to Partial Differential
Equations}. Berlin-Heidelberg-New York: Springer 1983.
%%%%%%%%%%%%%%%%%%%%%%%%%%%%%%%%%%%%%%%%%%%%%%%%%%%%%%%%%%%%
\bibitem{Okazawa73} N. Okazawa, \textit{Perturbations of Linear m-Accretive Operators}, Proc. Amer. Math. Soc. Vol. 37, No. 1 (Jan., 1973), pp. 169-174.
%%%%%%%%%%%%%%%%%%%%%%%%%%%%%%%%%%%%%%%%%%%%%%%%%%%%
%\bibitem{Okazawa69} N. Okazawa,\textit{ Two perturbation theorems for contraction semigroups in a Hilbert space}, Proc. Japan Acad. 45 (1969), 850-853.
%%%%%%%%%%%%%%%%%%%%%%%%%%%%%%%%%%%%%%%%%%%%%%%%%%%%%
%\bibitem{Okazawa77} N. Okazawa, \textit{Approximation of linear m-accretive operators in a Hilbert space},
%Osaka J. Math., 14 (1977), 85--94.
%%%%%%%%%%%%%%%%%%%%%%%%%%%%%%%%%%%%%%%%%%%%%%%%%%%%%%%%%
\bibitem{Okazawa82} N.Okazawa, \textit{On the perturbation of linear operators in Banach and Hilbert spaces}, J. Math. Soc. Japan 34 (1982) 677--701.
%%%%%%%%%%%%%%%%%%%%%%%%%%%%%%%%%%%%%%%%%%%%%%%%%%%%%%%%%%%
%\bibitem{Okazawa2002}N. Okazawa, \textit{Perturbation theory for m-accretive operators and
%generalized complex Ginzburg-Landau equations}, J. Math. Soc. Japan Vol. 54 No. 1 (2002), 1--19.
%%%%%%%%%%%%%%%%%%%%%%%%%%%%%%%%%%%%%%%%%%%%%%%%%%%%%%%%%%%%
\bibitem{Ota1984} S. \^Ota, \textit{Closed linear operators with domain containing their range}, Proc. Edinburgh Math. Soc. 27(1984), 229--233.
%%%%%%%%%%%%%%%%%%%%%%%%%%%%%%%%%%%%%%%%%%%%%%%%%%%%%%%%%%%%%%%%%%%%%
\bibitem{Shkalikov89} A. A. Shkalikov, \textit{Strongly damped pencils of operators and solvability of the corresponding operator-differential equations},  Math. USSR Sb. 63 (1989), 97-119.
%%%%%%%%%%%%%%%%%%%%%%%%%%%%%%%%%%%%%%%%%%%%%%%%%%%%%
%\bibitem{Sobajima14} M. Sobajima , \textit{A class of relatively bounded perturbations for generators of bounded analytic semigroups in Banach spaces }, J. Math. Anal. Appl. 416 (2014) 855--861
%%%%%%%%%%%%%%%%%%%%%%%%%%%%%%%%%%%%%%%%%%%%%%%%%%%%%
%\bibitem{Sohr81} H. Sohr,  \textit{Ein neues Surjektivitatskriterium im Hilbertraum}. Mh. Math. 91, 313--337
%(1981). 
%%%%%%%%%%%%%%%%%%%%%%%%%%%%%%%%%%%%%%%%%%%%%%%%%%%%%%%%
\bibitem{Stampfli67} J. G. Stampfli, \textit{Minimal range theorems for operators with thin spectra}, Pac. J. Math., 23 (1967), 601--612.
%%%%%%%%%%%%%%%%%%%%%%%%%%%%%%%%%%%%%%%%%%%%%%%%%%%
%\bibitem{Wust1971} R. Wust, \textit{Generalisations of Rellich's theorem on perturbation of (essentially) selfadjoint
%operators}, Math. Z. 119 (1971), 276--280.
%%%%%%%%%%%%%%%%%%%%%%%%%%%%%%%%%%%%%%%%%%%%%%%%%%%%%
\bibitem{Yoshikawa72} A. Yoshikawa, \textit{On Perturbation of closed operators in a Banach space}, J. Fac. Sci. Hokkaido Univ., 22 (1972), 50--61.
%%%%%%%%%%%%%%%%%%%%%%%%%%%%%%%%%%%%%%%%%%%%%%%%%
\bibitem{Yosida65} K. Yosida,\textit{ A perturbation theorem for semigroups of linear operators}, Proc. Japan
Acad. 41 (1965), 645--64.
%%%%%%%%%%%%%%%%%%%%%%%%%%%%%%%%%%%%%%%%%%%%%%%%%%%%%%%%%
\end{thebibliography}
\end{document}